\documentclass[a4paper,12pt]{article}
\usepackage{graphicx}
\usepackage[mathscr]{eucal}
\usepackage{multirow}
\usepackage{amsmath,amssymb,amsthm,cite,caption}
\usepackage{enumitem}
\usepackage{color}
\usepackage[labelformat=simple]{subcaption}

\newtheorem{theorem}{Theorem}[section]
\newtheorem{proposition}[theorem]{Proposition}
\newtheorem{lemma}[theorem]{Lemma}

\newtheorem{hypothesis}[theorem]{Hypothesis}

\def\real{\mathbb{R}}
\def\nat{\mathbb{N}}

\def\vp{\varphi}

\def\del{\nabla}

\def\play{\mathfrak{p}}
\def\stop{\mathfrak{s}}

\def\nato{\mathbb{N}\cup \{0\}}
\def\SS{\mathcal{S}}

\def\bbv{v}

\def\bbx{x}
\def\bbc{\beta}
\def\bbxi{p}

\def\bbphi{\varphi}
\def\FF{\mathcal{F}}
\def\oi{^{(i)}}
\def\oj{^{(j)}}
\def\on{^{(n)}}

\def\for{\mbox{ for }\ }
\def\pif{\mbox{ if }\ }

\newfont{\ctv}{msam10}

\newcommand{\bbox}{\mbox{\ctv \symbol{4}}}

\def\QED{{$\hfill\bbox$}}

\newenvironment{pf}[1]{\par\vskip1mm{\noindent\it #1.}\ }{\QED\par\vskip2mm}

\def\be{\begin{equation}\label}
\def\ee{\end{equation}}

\def\barr{\begin{array}}
	\def\earr{\end{array}}

\def\ber{\begin{eqnarray}}
\def\eer{\end{eqnarray}}

\def\bers{\begin{eqnarray*}}
	\def\eers{\end{eqnarray*}}

\def\bpf{\begin{pf}}
	\def\epf{\end{pf}}

\title{Global stability of a piecewise linear macroeconomic model with a continuum of equilibrium states and sticky expectation}

\author{
Pavel Krej\v{c}\'\i\thanks{Institute of Mathematics of the Czech Academy of Sciences, Czech Republic} \and
	Harbir Lamba\thanks{Department of Mathematical Sciences,
    George Mason University, VA, USA} \and Dmitrii~Rachinskii\thanks{Department of Mathematical Sciences, The University of Texas at Dallas, TX, USA}}

\date{}

\begin{document}

\rmfamily

\maketitle

\begin{abstract}
	We consider piecewise linear discrete time macroeconomic models, which possess a continuum of equilibrium states.
	These systems are obtained by replacing rational inflation expectations with a boundedly
	rational, and genuinely sticky, response of agents to changes in the actual inflation rate in a standard 
	Dynamic Stochastic General Equilibrium model. Both for a low-dimensional variant of the model, with one representative agent, and the multi-agent model, we show that, when exogenous
	noise is absent from the system,  the continuum of equilibrium states is the global attractor.
	Further, when a uniformly bounded noise is present, or the equilibrium states are destabilized by an imperfect Central Bank policy (or both), we estimate the size of the domain that attracts all the trajectories.
	The proofs are based on introducing a family of Lyapunov functions and, for the multi-agent model, deriving a formula for the inverse of the Prandtl-Ishlinskii operator acting in the space of discrete time inputs and outputs.
\end{abstract}


\smallskip


\section{Introduction} \label{intro}

Notions of {friction} and {stickiness} are widely accepted to exist in organizations, economies, and financial systems.
Drawing intuition from mathematical physics, one can propose that a system with (dry) friction should have a {\em continuum} of equilibrium states. For example, an object can achieve an equilibrium on a curved surface at any point where the slope does not exceed the dry friction coefficient because friction balances the gravity. 
Notably,
various hard-to-explain empirical regularities found in micro, macro, and organizational economic data, such as { path-dependence, 
	permanence, hysteresis,  boom blessings and recession curses} 
can be accounted for by the presence of many meta-stable states
with the associated long timescale dynamics, which frictions introduce into an economic system.
In other words, frictional effects at the micro-level might aggregate to macro-level long term memory effects.

Based on this premise, we propose to test how the introduction of internal frictions 
 affects dynamics of well-established single equilibrium economic models. It is a fundamental question what forms of friction actually arise from the behavior of economic agents, and how these can be modeled.
In this paper, we take a phenomenological approach to modeling frictions using {\em play operators}, which 
are common in physics and engineering applications. In the context of economics, they are associated with the notions of `{threshold}' and `{inaction band}'. These operators have been described in several equivalent ways including piecewise linear (PWL) functions and variational inequalities. 
We proceed by positing that inflation and inflation expectation  are related by the play operator (or a combination thereof in the multi-agent variant of the model) in a version of a Dynamic Stochastic General Equilibrium (DSGE)
macroeconomic model. Our analysis focuses on one particular feature of the 
adapted model with frictions. Namely, we consider local and global stability of this discrete time PWL system, which naturally possesses a continuum of equilibrium states.

In the remainder of this section some economic background and the structure of the paper are briefly discussed.

\smallskip
{\bf Economic background and motivation.}
Seminal work by the likes of Walras and Jevons in the late 1800's steered the field of economics away from the domain of philosophy, and laid  foundations for it to develop into a mathematically tractable, exact science. Early pioneers of what became known as \textit{Neoclassical Economics} saw parallels between economic systems and the equilibrating forces in nature, and thus borrowed heavily from mathematical physics to derive both intuition and methodologies\cite{Mirowski}.  Despite concerns about the simplifying assumptions needed for this t\^{a}tonnement-driven view (e.g. \cite{Marshall,Keynes}), the single-equilibrium approach is still considered good enough in many settings, and has not been 
supplanted by any widely-accepted, complete system. 
Significant developments have come from introducing concepts like {Rational Expectations} and {Sticky Models}.
The former is the assumption of aggregate consistency in dynamic models \cite{m61}. This view admits that agent's expectations about future uncertainties may be wrong individually, but in aggregate are in agreement with the model itself. In other words, although the future is not fully predictable, agents' expectations are assumed not to be systematically biased and collectively use all relevant information in forming expectations of economic variables. Meanwhile the later includes the widely-used sticky models of Calvo \cite{calvo1983staggered} and the sticky-information of Mankiw and Reis \cite{mankiw2002sticky}. These models are concerned with similar observations about the way in which real life agents do not instantaneously move to the `correct' price or opinion but rather do so at a fixed rate and  can be represented mathematically by introducing a delay term into the  relevant equations. In the absence of noise the same optimal equilibrium solution will be reached as if the stickiness were absent. Continua of possible equilibria can also occur in such models (see for  example \cite{benhabib1999indeterminacy,evans2015observability}) but only in certain special cases (such as a passive interest-rate  policy \cite{calvo1983staggered,antinolfi2007monetary}) and are considered an extreme form of {\em indeterminacy}.

However, a robust empirical feature of economic output for developed countries is that output and output growth are non-normally distributed, exhibiting fat tails and excess kurtosis. In many models such booms and busts are explained by the occurrence of large (unpredictable) exogenous shocks followed by tranquil periods when nothing leads to non-normality.   Explaining such regularities within linear unique-equilibrium economic models often involves adding \textit{a posteriori} assumptions such as the existence of an eigenvalue with largest modulus close to, but inside, the unit circle (the Unit Root Hypothesis). A better model should generate non-normality of the output data from within the theory.  

The critique of unique-equilibrium models has a long history which we shall not attempt to detail here. For example, many have eloquently pointed out fundamental issues with the assumed equilibrating processes and the ways in which the ``aggregation problem" was being solved  (e.g. \cite{robinson1974,n72,s97a,colander2008beyond}).  This approach is taken by multi-agent models including models of agents' behavior under imperfect information, cognitive limitations, or endogenously generated ``animal spirits'' \cite{grauwe}. 
 Here, the themes of confusion, rational inattention, simplification, and bounded rationality take center stage (e.g. the work of Sims and Gabaix). Non-normality 
 in these models can be associated with cascading effects. Most of these models are however purely numerical, and incorporate assumptions inspired by the social sciences (e.g. psychology and sociology) in order to add more realistic behavior to the modeling construct of individual actors
 (much recent work has focused on the question of whether individual actors can ever be sufficiently rational).   

Our approach is complementary to many of the above-referenced views. 
In this paper, we argue that a lingering constraint 
of single-equillibrium in all these models is unnecessary.
 Hysteresis is a well-understood property that explains the ``stickiness" observed in many physical systems (e.g. plasticity, magnetism), and which explains how many equillibria can arise and be stable. Hysteresis has only recently been explicitly considered in economic time series such as the unemployment rate \cite{cross, cross0, scihys}.  The form of stickiness that we use is, to our knowledge, new in a economic setting and differs from, for example, the stickiness of the Calvo pricing model \cite{calvo1983staggered} where hypothetical agents are only allowed  to adjust (to the correct price) at a fixed rate.  The way in which we incorporate stickiness into the model will be justified and described more fully below but, briefly, our sticky variables can only be in one of two modes. They are either currently `stuck' at some value or they are being `dragged' along by  some other (related) variable because the maximum allowable difference between them has been reached. Hence,  our agents are truly stuck (not just delayed) until forced to adjust by the discrepancy with the actual inflation rate.  If an equilibrium is reached it is chosen by the prior states of the system, and a continuum of equilibria is an intrinsic feature of the model. 

The research into how expectations are formed is extensive but far
from conclusive, see for example
\cite{curtin2010inflation,rudd2006can,branch2007sticky,carroll2003macroeconomic,mankiw2003disagreement}. However
the idea of threshold effects and a `harmless interval' of inflation
is not new in economics
\cite{Threshold2010,kremer2013inflation,2013inflation,F2010inflation,2001threshold}.
A person may concurrently be subject to many or all of the limitations (e.g. rational inattention and bounded rationality and confusion, etc.), but for our purposes it is enough to assume that individual actors behave according to the play operator we describe, and we can remain relatively agnostic as to which mechanism may be driving the general features of the play operator (band of inaction, thresholds).  
In the absence of any exogenous forcing it would be very easy to
distinguish between Calvo-type stickiness and the stuck-then-dragged
behavior we investigate here --- indeed Calvo stickiness would most
likely be observed since agents could tell far more easily over time
that, for example, their wage demands were too low and they were
losing purchasing power. However, given the uncertainty of reality and
the very limited cognitive skills or interest in forecasting
of most economic agents, that may no longer hold.

\smallskip
{\bf Structure of the paper.}
In the next section, we present a discrete time 3-dimensional PWL macroeconomic model with sticky inflation expectations modeled by a play operator. The model uses the notion of a representative agent. When the exogenous noise terms are absent from the system, it has a line segment of equilibrium points. Section 3 contains main results. In particular, in the system without exogenous noise, a simple condition ensures that the line segment of equilibrium points is the global attractor (Section 3.1). In the presence of uniformly bounded noise, we obtain an estimate of the globally attracting domain, which is proportional to the supremum norm of the noise. Interestingly, this estimate is uniform with respect to the parameter that controls the amount of stickiness in the expectation of future inflation rate. 
We then consider further variants of the model. First, we add stickiness
into the response of the Central Bank (Section 3.2). This can destabilize the equilibrium states (leading to periodic, quasiperiodic or more complex dynamics \cite{6authors} with the associated border collision bifurcations, which are typical of piecewise smooth systems \cite{Gallery}) but the system still possesses a bounded globally attracting domain. 
Then we consider a multi-agent variant of the model (Section 3.3). 
This is an $(n+2)$-dimensional PWL system with $2n$ switching surfaces.
We show that the $n$-dimensional continuum of equilibrium states of this system is the global attractor.
Finally, the proofs based on constructing a family of Lyapunov functions are presented in Section 4.
In order to apply the Lyapunov function to the multi-agent model, we adapt 
 a technique from the theory of hysteresis operators \cite{KP, BS}.
 Namely, an explicit formula for the inverse of the Prandtl-Ishlinskii operator acting in the space of discrete time inputs and outputs is 
 derived and used. 
We conclude with a summary of the main results and some suggestions for future work.

\section{The model} \label{mode}

\subsection{DSGE modeling framework} 
The standard approach to the problem of
aggregating expectations is to introduce a `Representative Agent'
whose expectations are fully-informed and rational and consistent with
the model itself. Here, an aggregation of {\em boundedly}
rational agents into a similar Representative  is required.
Our approach is similar in spirit to that of De Grauwe
\cite{grauwe}
but we use a different model of boundedly rational agents' behavior.

We start from a dynamic stochastic general equilibrium (DSGE)
macroeconomics model, which includes aggregate demand and aggregate
supply equations
\begin{equation}\label{e1}
\begin{array}{l}
x_t=b_1p_t+(1-b_1)x_{t-1}+b_2y_t+\eta_t,\\
 y_t=(1-a_1)y_{t-1}-a_2(r_t-p_t)+\epsilon_t,\\
\end{array}
\end{equation}
augmented with the interest rate-setting Taylor rule 
\begin{equation}\label{eqn:M1'}
r_t=c_1x_t+c_2y_t+\xi_t,
\end{equation}
where $y_t$ is output gap (or unemployment rate, or another measure of
economic activity such as gross domestic product), $x_t$ is inflation
rate, $r_t$ is interest rate, $p_t$ is the economic agents' aggregate
expectation of future inflation rate, $\eta_t$, $\epsilon_t$, $\xi_t$ are
exogenous noise terms, and $t \in \nat$. The parameters satisfy 
\[
0\le a_1 <1, \qquad 0<b_1<1, \qquad a_2, b_2, c_1, c_2 > 0.
\]
Of specific interest is the case $a_1=0$.

This model is close to the 
model used in \cite{grauwe} but simpler in that we do not
include the aggregate expectation of the output gap and the
correlation between the subsequent values of the interest rate.
The inclusion of such factors does not affect our most
significant qualitative observations, but would complicate some
aspects of the rigorous analysis that we present.\footnote{More complicated variants of DSGE models, employing many more variables to represent different sectors of the economy, are widely used by central banks to help determine interest rate policy \cite{dsge6,dsge4,dsge1,dsge2}. DYNARE software platform supporting both DSGE models relying on the rational expectations hypothesis and some models where agents have limited rationality or imperfect knowledge of the state of the economy is available at \cite{dynare}.} 

The novelty of our modeling strategy is in how we define the
relationship between the aggregate expectation of inflation $p_t$ and
the inflation rate $x_t$.  

\subsection{Sticky expectation of inflation}\label{sec_play}
We start from the empirical evidence cited above that individual agents'
expectations are often sticky and may lag behind the currently
observable values before they start to move. We also posit that
this gap between future expectations and current reality cannot grow
too large. We then imbue our now boundedly rational Representative
Agent with these same properties.
More precisely, 
we assume the following rules that define the variations of the
expectation of future inflation rate $p_t$ with the actual inflation
rate $x_t$ at integer times $t$:
\begin{itemize}
\item[(i)] The value of the difference $|p_t-x_t|$ never exceeds a certain bound $\rho$;

\item[(ii)] As long as the above restriction is satisfied, the expectation does not change, i.e. $|x_t-p_{t-1}|\le \rho$ implies $p_t=p_{t-1}$;

\item[(iii)] If the expectation has to change, it makes the minimal
  increment consistent with constraint (i).
\end{itemize}

Rule (ii) introduces stickiness in the dependence of $p_t$ on $x_t$,
while (i) states that the expected inflation rate cannot
deviate from the actual rate more than prescribed by a threshold value
$\rho$.  Hence $p_t$ follows $x_t$ reasonably closely but on the other
hand is conservative because it remains indifferent to variations of
$x_t$ limited to a (moving) window $p-\rho\le x\le p+\rho$.  The last
rule (iii) enforces continuity of the relationship between $p_t$ and
$x_t$ and, in this sense, can be considered as a technical modeling
assumption that is mathematically convenient.

Rules (i)--(iii) are expressed by the formula
\begin{equation}\label{formula}
p_{t}=x_t+\Phi_\rho(p_{t-1}-x_t),\qquad t\in\mathbb{N}
\end{equation}
with the piecewise linear saturation function 
\begin{equation}\label{formula'}
\Phi_\rho(v)=\left\{\begin{array}{rlc} 
\rho & {\rm if} & v\ge \rho,\\
v & {\rm if} & -\rho<v<\rho,\\
-\rho & {\rm if} & v\le -\rho.
\end{array}
\right.
\end{equation}
Equations \eqref{e1} and \eqref{eqn:M1'}, completed with formulas
\eqref{formula} and \eqref{formula'}, form a closed 3-dimensional PWL model for the
evolution of the aggregated variables $x_t,y_t, p_t$.
Another variant of this model has been considered in \cite{6authors}.
Iterations of system \eqref{e1}--\eqref{formula'} start
from a set of initial conditions $x_0,y_0,p_0$ satisfying 
$|x_0-p_0|\le \rho$.

Some further terminology will be useful. Denote by $\SS$ the set of all real sequences $\bbx = (x_0, x_1, \dots)$. For a given parameter $\rho>0$, the {\em play operator\/} $\play_\rho: [-\rho,\rho]\times \SS \to \SS$ is defined
as the mapping which with a given initial condition $s_0 \in [-\rho, \rho]$ and a sequence $\bbx \in \SS$
associates the sequence 
\begin{equation}\label{4}
p=(p_0,p_1,\ldots)=\play_{\rho}[s_0,x]_t\in \SS
\end{equation} 
according to the formula \eqref{formula} with $p_0=x_0-s_0$ \cite{KP}.
The parameter $\rho$ is called the {\em threshold}, see Fig.~1 (left).
A dual mapping  $\stop_\rho: [-\rho,\rho]\times \SS \to \SS$, which is defined by
the relationship 
\begin{equation}\label{formula''}
s_{t}=\Phi_\rho(x_t-x_{t-1}+s_{t-1}),\qquad t\in\mathbb{N}
\end{equation}
for an arbitrary pair $(s_0,x)\in [-\rho,\rho]\times \SS$, is known as the {\em stop operator}
\begin{equation}\label{4'}
s=(s_0,s_1,\ldots)=\stop_{\rho}[s_0,x]_t\in \SS,
\end{equation} 
see Fig.~1 (right). By definition the play and stop operators sum up to the identity:
\[
p_t+s_t=x_t
\]
for the sequences \eqref{4}, \eqref{4'}.

\begin{figure}
	\centering
	\includegraphics[width=0.45\textwidth]{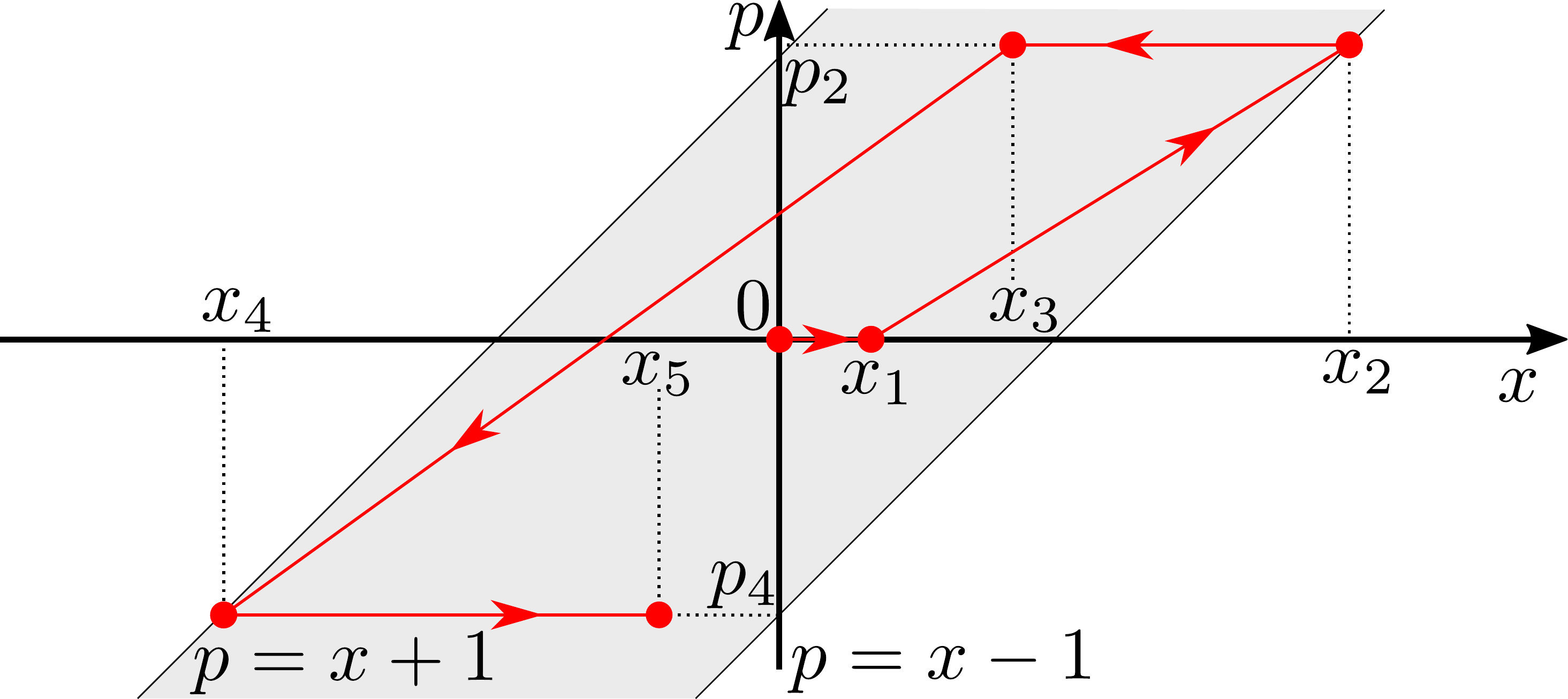}
	\qquad
	\includegraphics[width=0.45\textwidth]{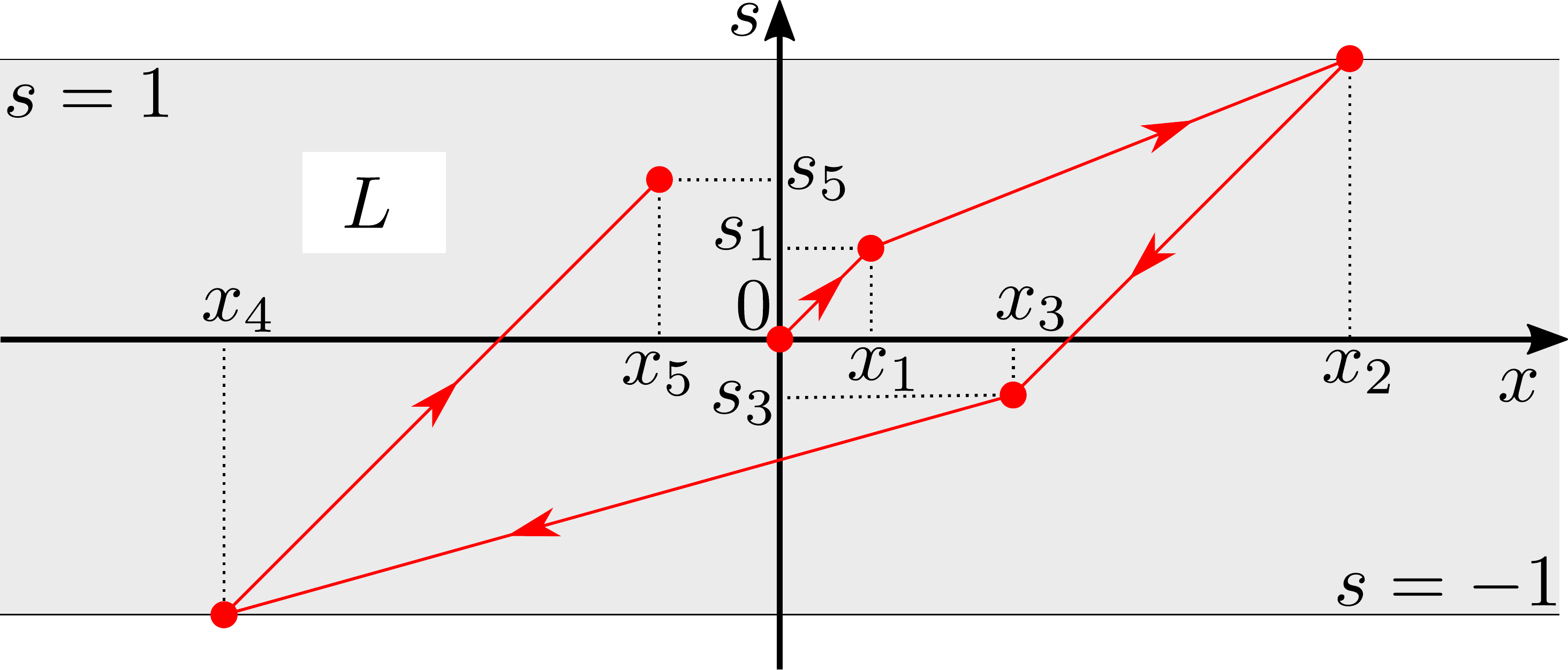}
	\caption{\small 
		Graphical representation
		of the play and stop operators with the threshold $\rho=1$. Left: Example of an input-output sequence $(x_t,p_t)$ of the play
		operator. For each $t$,\ $(x_t,p_t)$ is the nearest point to $(x_{t-1},p_{t-1})$, which belongs to the band $|x-p|\le 1$ (shown in grey) and has $x=x_t$. Right: The corresponding input-output sequence $(x_t,s_t)=(x_t, x_t-p_t)$ of the stop operator satisfies $|s_t|\le 1$ for all $t$.}
	\label{fig1}
\end{figure}

\begin{figure}
	\centering
	\includegraphics[width=0.40\textwidth]{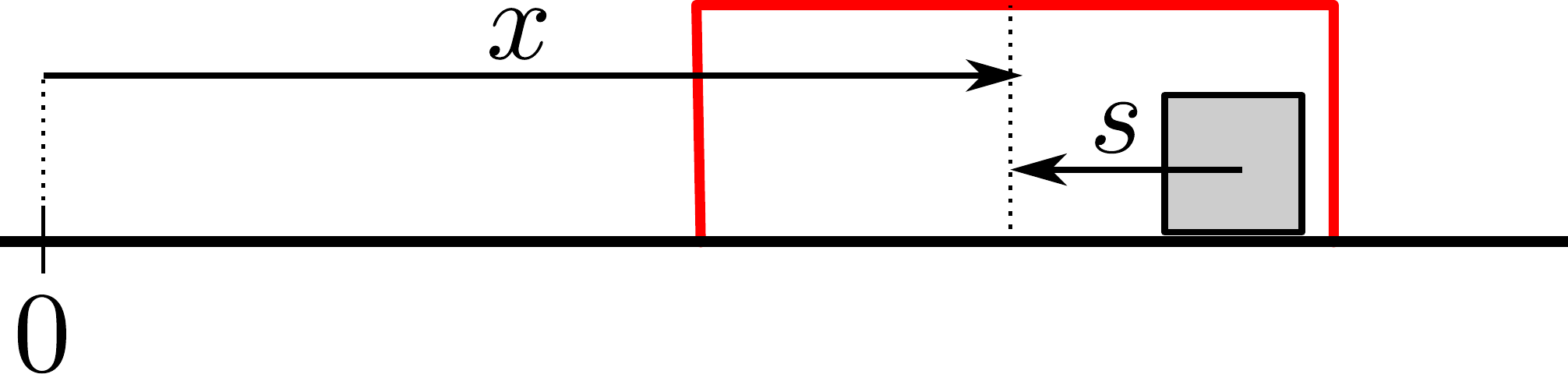}\qquad\quad
	\includegraphics[width=0.40\textwidth]{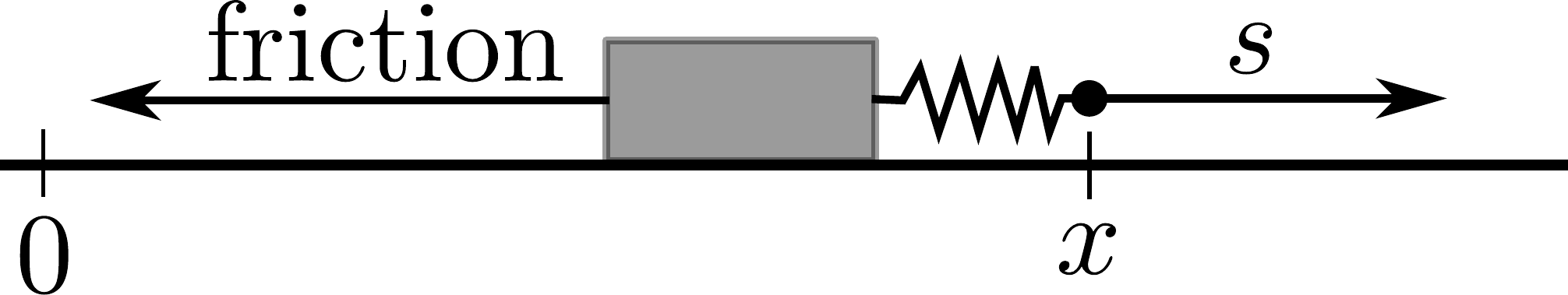}
	\caption{\small Continuous time mechanical interpretations of the stop operator. Left: 
		The position $x$ of the center of the frame, which has length 2$\rho$,  is the input; the relative
		position $s$ of the center of the frame  with respect to
		the center of the grey box is the output. The frame moves and drags the
		box creating the input-output sequence $(x,s)$.
		Right: Schematic of the
		Prandtl's model of quasi-static elastoplasticity \cite{prandtl}. The grey box is not moving
		unless
		the absolute value of the force $s$ of the ideal spring reaches the
		maximal value $\rho$ of the Coulomb friction force between the box and the surface (equivalently, $s$ can denote the elongation of the spring); $x$ is the coordinate of the end of the spring.
	}
	\label{fig:test}
	\vskip-4truemm
\end{figure}

One can think of the play operator as having two modes, see Fig.~2 (left). A `stuck mode'
where it will not respond to small changes in the input $x_t$ and a
`dragged mode' where the absolute difference between the input $x_t$ and output $p_t$ are
at the maximum allowable and changes to the input, in the correct
direction, will drag the output along with it.
Further, in the context of our model, the output of the stop operator $s_t$
measures the difference between the inflation rate and the expectation
of the future inflation rate, hence $s_t$ remains within the bound
$|s_t|\le \rho$ at all times. We shall refer to the variable $s_t = x_t-p_t$ as the {\em perception gap}. Interestingly the explicit relationship
\eqref{formula''} has been observed in actual economic data \cite{g, laura}.
Fig. 2 (right) gives an interpretation of the stop operator combing an ideal spring and a dry friction element as used in mechanics.

\section{Main results}

\subsection{Autonomous system}

We first consider system \eqref{e1}--\eqref{formula} without noise terms:
\begin{equation}\label{e1'}
\begin{array}{l}
x_t=b_1p_t+(1-b_1)x_{t-1}+b_2y_t,\\
 y_t=(1-a_1)y_{t-1}-a_2(c_1x_t+c_2y_t-p_t),\\
p_{t}=x_t+\Phi_\rho(p_{t-1}-x_t).
\end{array}
\end{equation}
Equilibrium points of this system form a line segment
\be{aA}
{\mathfrak A}=\left\{(x_*,y_*,p_*)(u):\ x_*(u)=\kappa u,\ y_*(u)=\frac{b_1}{b_2}u,\ p_*(u)=(\kappa-1)u,\ |u|\le \rho\right\}
\ee
with
\[
\kappa=\frac{a_2(b_2+b_1c_2)+a_1b_1}{a_2b_2(1-c_1)}.
\]

\begin{theorem}
	\label{t1}
	If $c_1>1$, then the line segment $\mathfrak A$ of equilibrium points is the global attractor for system \eqref{e1'}.
	Further, any trajectory converges to an equilibrium point $(x_*,y_*,p_*)\in \mathfrak A$.
\end{theorem}

We note that system \eqref{e1'} is written in an implicit form. If $c_1>1$ as in Theorem \ref{t1},
it is easy to rewrite this system as an explicit PWL map $(x_t,y_t,p_t)=f(x_{t-1},y_{t-1},p_{t-1})$.

\subsection{Sticky Central Bank response}
The Central Bank policy can presumably exhibit stickiness too. To explore this scenario
in this section we replace the Taylor rule \eqref{eqn:M1'} with the relation
\be{de1}
r_t= \play_\sigma [c_1 x + c_2 y]_t + \xi_t
\ee
also involving a play operator
with an initial condition $r_0$ such that $|r_0 - (c_1 x_0 + c_2 y_0)| \le \sigma$.
The play operator $\play_\sigma$ with a threshold $\sigma\ge 0$ independent of $\rho$ in \eqref{de1}
should express the fact that the central bank's decisions do not immediately follow
the instantaneous value of $c_1 x_t + c_2 y_t$, but they are activated only if the difference
between $r_t$ and $c_1 x_t + c_2 y_t$ risks to exceed a given value $\sigma$.
For $\sigma = 0$, $\play_\sigma$ is the identity mapping and 
\eqref{de1} becomes \eqref{eqn:M1'}.
It is worth noting that
for $\sigma > 0$ and a sequence $\{v_t\}$, $r_t = \play_\sigma[v]_t$ is the sequence with minimal
variation in the $\sigma$-neighborhood of $\{v_t\}$, that is, the implication
\be{mvar}
\hat r_0 = r_0, \ |\hat r_t - v_t| \le \sigma \ \ \forall t \in \nat \cup \{0\}
\quad \Rightarrow \quad \sum_{t=1}^T |\hat r_t - \hat r_{t-1}|
\ge \sum_{t=1}^T |r_t - r_{t-1}| 
\ee
holds for every sequence $\{\hat r_t\}$ and any $T \in \nat$ as a special case of \cite[Theorem 1.2 and Corollary 1.6]{klf}.
The meaning of \eqref{de1} can thus be interpreted as the optimal
strategy allowing to stay in a $\sigma$-neighborhood of the input sequence under minimal ``cost'' incurred by variations of the interest rate.

Re-writing equation \eqref{de10} equivalently as
\begin{equation}\label{formulaa}
r_{t}=c_1 x_t + c_2 y_t+\Phi_\sigma(r_{t-1}-c_1 x_{t-1} - c_2 y_{t-1})+\xi_t,
\end{equation}
we obtain a 4-dimensional PWL system \eqref{e1}, \eqref{formula}, \eqref{formulaa}.
In the absence of noise, i.e. when $\eta_t=\epsilon_t=\xi_t=0$ for all $t$, this system has the following set of equilibrium points:
\[
{\mathfrak M}=\left\{(x_*,y_*,p_*,r_*):\ x_*=\kappa u+\frac{v}{c_1-1},\ y_*=\frac{b_1}{b_2}u,\ p_*=(\kappa-1)u+\frac{v}{c_1-1},\right.
\]
\[
\left.  r_*=\left(c_1\kappa + \frac{b_1c_2}{b_2}\right)u+\frac{v}{c_1-1},\quad |u|\le \rho, \ |v|\le \sigma\right\}.
\]
The subset 
\be{a0}
{\mathfrak A_0}=\left\{(x_*,y_*,p_*,r_*)\in {\mathfrak M}: \ v=0
\right\}
\ee
of $\mathfrak M$ can be considered as a natural embedding of the equilibrium set \eqref{aA}
into the 4-dimensional phase space of system \eqref{e1}, \eqref{formula}, \eqref{formulaa}.

It should be pointed out right away that stickiness in the Taylor rule can have a destabilizing effect on the
equilibrium states. To see this, consider for simplicity the case when $\rho=0$, which implies $p=x$ and $s=0$, i.e. we remove stickiness in the inflation expectation. Note that the system is locally linear in a vicinity of every equilibrium point belonging to $\mathfrak M$
and satisfying $ |v| < \sigma$. Further, a simple calculation 
shows that the determinant of the linearization at such equilibrium points
equals $(1-a_1)(1-b_1)/(1-b_1-a_2b_2)$. 
If $1-b_1>a_2 b_2>a_1(1-b_1)$, then this determinant is greater than $1$, hence these equilibria are unstable (in particular, they are unstable in the important case of $a_1=0$, $1-b_1>a_2 b_2$).
Numerical examples of several attractors including periodic orbits of different periods, a quasiperiodic orbit or a union of two equilibrium points corresponding to $v=\pm\sigma$ (end points of the line segment of equilibrium points) can be found in \cite{}.

The goal of this section is to estimate how far trajectories 
can deviate from the equlibrium set $\mathfrak A_0$ due to stickiness in the Taylor rule and a uniformly bounded noise.

\begin{theorem}\label{t2}
	Let us consider system \eqref{e1}, \eqref{formula}, \eqref{formulaa} with uniformly bounded exogenous terms $\eta_t,\epsilon_t,\xi_t$ and with $c_1>1$.
	There exist constants $L_1,L_2$, which depend on the parameters $a_{1},a_{2}, b_{1},b_{2}, c_{1},c_{2}$ but are independent of the threshold $\rho$ of the inflation expectation (see \eqref{formula}), such that every trajectory satisfies
	\be{estdi}
	\limsup_{t\to\infty} |x_t-x_*(s_t)|, \
	\limsup_{t\to\infty} |y_t-y_*(s_t)|\le L_1 \sigma + L_2
m,	\ee
where $(x_*,y_*,p_*)(\cdot)$ is defined in \eqref{aA}, $s_t=x_t-p_t$, $\sigma$ is the threshold of the play operator in the central bank's policy (see \eqref{de1})  and
\begin{equation}\label{added}
m=\sup_t \max\{|\epsilon_t|, |\eta_t|, |\xi_t|\} <\infty.
\end{equation}
	
\end{theorem}

The proof presented below suggests an explicit upper bound for the coefficients $L_1, L_2$ in \eqref{estdi}.
Due to relations $p=x-s$ and \eqref{formulaa}, estimates \eqref{estdi} imply
similar estimates for $p_t$ and $r_t$: 
\be{es1di}
	\limsup_{t\to\infty} |p_t-p_*(s_t)|
\le L_1 \sigma + L_2m,
\ee
\be{es2di}
\limsup_{t\to\infty} |r_t-c_1x_*(s_t)-c_2y_*(s_t)|\le (c_1+c_2)(L_1\sigma+L_2m)+\sigma+m.
\ee

According to Theorem \ref{t2}, estimates \eqref{estdi}--\eqref{es2di} are uniform with respect to $\rho$.
In particular, if $\sigma=0$ and hence system \eqref{e1}, \eqref{formula}, \eqref{formulaa} becomes equivalent to system \eqref{e1}--\eqref{formula}, then every trajectory converges to a neigborhood 
\[
{\mathfrak A}(R)=\left\{(x,y,p): \ \min_{(x_*,y_*,p_*)\in{\mathfrak A}} \max\{|x-x_*|,|y-y_*|, |p-p_*|\}\le R\right\}
\]
of the set \eqref{aA} of equilibrium points, and the size of this neighborhood is proportional to the supremum norm of the noise terms, $R=L_2 m$, and is independent of the threshold $\rho$ of inflation expectations (cf. Theorem \ref{t1}). 

\subsection{A multi-agent model}\label{multi}
Model \eqref{e1}--\eqref{formula}
can be easily extended to account for 
differing types of agent with different inflation rate expectation thresholds.
To this end, we replace the simple relationship \eqref{4}
between $p_t$ and $x_t$ (which is equivalent to \eqref{formula}) with the equation
\be{m1}
p_t = \sum_{i=1}^n \nu_i \play_{\rho_i}[\beta_i, x]_t.
\ee
Here the play operator $\play_{\rho_i}$ models the expectation
of inflation by the $i$-th agent; $p_t$ is the aggregate expectation of inflation; 
$\nu_i>0$ is a weight measuring the contribution of agent's expectation 
of inflation to the aggregate quantity; $\rho_i$ is an individual threshold characterizing
the behavior of the $i$-th agent; $\beta_i$ is the initial condition for the corresponding play operator; and we assume the ordering $0 < \rho_1 < \dots < \rho_n$. 

Operator \eqref{m1} is known as a (discrete) Prandtl-Ishlinskii (PI) operator  \cite{AYIshlinskiiStatMethods,prandtl,KP}.\footnote{The continuous time counterpart of this
operator is widely used as a friction model in mechanical applications \cite{Lam} as well as for modeling elastoplastic systems \cite{pi0}, constitutive laws of smart materials \cite{pra}, and material fatigue \cite{Ryc}.} 
In analysis of this operator, it is convenient to 
restrict the set of initial conditions of the play operators in \eqref{m1}.
In what follows, we assume that 
$|\beta_1| < \rho_1$ and $|\beta_i - \beta_{i-1}| \le \rho_i - \rho_{i-1}$ for all $2 \le i \le n$.
Further, the coefficients $\nu_i$ are assumed to fulfill the condition
\be{m2}
\nu_i > 0 \ \for i=1, \dots, n\,, \quad \nu_0 := 1 -\sum_{i=1}^n \nu_i \ge 0\,.
\ee
It will be useful to define the quantity $s_t$ similarly as in \eqref{4'}, that is
\be{m3}
s_t = x_t - p_t = \nu_0 x_t + \sum_{i=1}^n \nu_i \stop_{\rho_i}[\beta_i, x]_t\,,
\ee
where $\stop_{\rho_i}[\beta_i, x]_t = x_t - \play_{\rho_i}[\beta_i, x]_t$ for all $t\in \nato$.

Equations \eqref{e1}, \eqref{eqn:M1'}, \eqref{m1} form an $(n+2)$-dimensional PWL system.
Similarly, if the Taylor rule \eqref{eqn:M1'} is replaced with its sticky counterpart \eqref{formulaa}, we obtain a PWL system of dimension $n+3$. One can formulate natural analogs of Theorems \ref{t1}, \ref{t2} for these systems. For example, let us consider the analog of Theorem \ref{t1} for the autonomous system
\begin{equation}\label{e1'''}
\begin{array}{l}
x_t=b_1p_t+(1-b_1)x_{t-1}+b_2y_t,\\
y_t=(1-a_1)y_{t-1}-a_2(c_1x_t+c_2y_t-p_t)
\end{array}
\end{equation}
coupled with formula \eqref{m1} for the aggregate expectation of inflation.

\begin{theorem}
	\label{t3}
	If $c_1>1$, then any trajectory of system \eqref{m1}, \eqref{e1'''} converges to an equilibrium point of this system.
\end{theorem}

We note that equilibrium points of system \eqref{m1}, \eqref{e1'''} form an $n$-dimensional parallelepiped in its phase space.
The proof of Theorem \ref{t3} uses the constructions from the proof of Theorem \ref{t1} but additionally relies on an inversion formula for the Prandtl-Ishlinskii operator. This proof also shows a possible way to extend Theorem \ref{t2} to the multi-agent model \eqref{e1}, \eqref{formulaa}, \eqref{m1} with sticky inflation expectation and exogenous noise.


\section{Proofs}

\subsection{Play and stop operators}\label{pla}

For the reader's convenience, we summarize here some well-known properties of the discrete time stop operator
$s_t = \stop_\rho[s_0,x]_t$ and play operator
$p_t = \play_\rho[s_0,x]_t$ which are needed in the sequel.

\begin{lemma}\label{l1}
	Let $\{x_t; t\in \nat\cup\{0\}\}$ be a given sequence. Then $p_t, s_t$ satisfy \eqref{4}, \eqref{4'} if and only if
	$|s_t| \le \rho$ for all $t\in \nat\cup\{0\}$ and the variational inequality
	\be{vari}
	(p_t - p_{t-1}, s_t - z) \ge 0 \quad \forall t \in \nat\cup\{0\}
	\ee 
	holds for every $z \in [-\rho, \rho]$.
\end{lemma}

\bpf{Proof}
Relations \eqref{4}, \eqref{4'} are equivalent to the series of implications
$$
p_t - p_{t-1} > 0 \ \Rightarrow \ p_t = x_t - \rho  \ \Rightarrow \ s_t = \rho\,, \
p_t - p_{t-1} < 0 \ \Rightarrow \ p_t = x_t + \rho  \ \Rightarrow \ s_t = -\rho\,,
$$
which is in turn equivalent to \eqref{vari} under the condition $|s_t| \le \rho$ for all $t\in \nat$.
\epf

For a generic sequence $\{z_t; t\in \nato\}$, we introduce the notation
\be{de2}
\del_t z := z_t - z_{t-1}\quad \pif t \ge 1; \quad \del_t^2 z = z_t - 2z_{t-1} + z_{t-2} \quad \pif t \ge 2\,.
\ee
We will systematically use the identity
\be{de7}
\del_t z\, z_t  = \frac12 \big(z_t^2 - z_{t-1}^2\big) + \frac12 (\del_t z)^2.
\ee
Choosing in \eqref{vari} the value $z = s_{t-1}$, we obtain that $\del_t p \del_t s \ge 0$, hence
\be{de4}
\del_t x \del_t s \ge (\del_t s)^2\,.
\ee
Furthermore, summing the inequalities
\bers
(p_t - p_{t-1})(s_t - s_{t-1}) &\ge& 0\,,\\
(p_{t-1} - p_{t-2})(s_{t-1} - s_t) &\ge& 0\,,
\eers
which follow from \eqref{vari} by the choice $z = s_{t-1}$ and $z = s_t$, respectively,
we obtain that $\del_t^2 p \del_t s \ge 0$, hence
\be{de5}
\del_t^2 x \del_t s \ge \del_t^2 s \del_t s = \frac12 \big((\del_t s)^2 - (\del_{t-1}s)^2\big)
+ \frac12 (\del_t^2 s)^2\,,
\ee
and similarly
\be{de6}
\del_t^2 x \del_t x  = \frac12 \big((\del_t x)^2 - (\del_{t-1}x)^2\big)
+ \frac12 (\del_t^2 x)^2\,,
\ee
which is a special case of identity \eqref{de7} with $z_t = \del_t s$ and $z_t = \del_t x$.

\begin{lemma}\label{l2}
	For a given sequence $\{x_t; t\in \nato\}$, put $p_t = \play_\rho[s_0,x]_t$, $s_t = x_t - p_t$ with some given initial condition
	$s_0\in [- \rho,  \rho]$.
	Let $q_t = x_t + \delta s_t = (1+\delta)x_t - \delta p_t$ for some $\delta > -1$. Then
	\be{inve}
	x_t = \frac{1}{1+\delta}\left(q_t + \delta \play_{(1+\delta)\rho}[(1+\delta)s_0,q]_t\right)\,.
	\ee
\end{lemma}

\bpf{Proof}
We have $q_t - p_t = (1+\delta)s_t$, hence $|q_t - p_t| \le (1+\delta)\rho$, and
\be{variq}
(p_t - p_{t-1},q_t - p_t  - (1+\delta)\rho x) \ge 0 \quad \forall t \in \nat \ \forall |x|\le 1\,.
\ee 
By Lemma \ref{l1}, this implies that $p_t = \play_{(1+\delta)\rho}[(1+\delta)s_0,q]_t$ and the assertion follows.
\epf

\begin{lemma}\label{l3}
	For a given sequence $\{x_t; t\in \nato\}$, put $p_t = \play_\rho[s_0,x]_t$ with some given initial condition
	$s_0\in [ - \rho, \rho]$. Then for every $t,j \in \nat$ we have
	\be{lip}
	|p_{t+j} - p_{t}| \le \max_{i=1, \dots, j}\{|x_{t+i} - x_t|\}\,.
	\ee 
\end{lemma}

\bpf{Proof}
We fix $t\in \nato, J \in \nat$ and for $j = 0,1,\dots,J$ set
$$
S_j = \max\{|p_{t+j} - p_{t}|^2, \max_{i=1, \dots, J}\{|x_{t+i} - x_t|^2\}\}\,.
$$
The proof will be complete if we prove that the sequence $\{S_j\}$ is nonincreasing
for $j = 0,1,\dots,J$. Indeed, then $S_J \le S_0$, which is precisely the desired statement.

Assume for contradiction that $S_j > S_{j-1}$ for some $j = 1,\dots,J$. Then
\ber \label{c1}
|p_{t+j} - p_{t}| &>& \max_{i=1, \dots, J}\{|x_{t+i} - x_t|\}\,,\\  \label{c2}
(p_{t+j} - p_{t})^2 &>& (p_{t+j-1} - p_{t})^2\,.
\eer
Inequality \eqref{c2} can be equivalently written as
\be{c3}
(p_{t+j} - p_{t+j-1})(p_{t+j} - p_{t}) > \frac12(p_{t+j} - p_{t+j-1})^2 > 0\,.
\ee
We now replace in \eqref{vari} written for $t+j$ instead of $t$ the element $z$ by $s_t$
and obtain
\be{c4}
(p_{t+j} - p_{t+j-1})(s_{t+j} - s_{t}) \ge 0\,,
\ee
hence, combining \eqref{c3} with \eqref{c4}, we have
$$
(p_{t+j} - p_{t})(s_{t+j} - s_{t}) \ge 0\,,
$$
which implies that
$$
(p_{t+j} - p_{t})^2 \le (p_{t+j} - p_{t})(x_{t+j} - x_{t})
$$
in contradiction with \eqref{c1}. This completes the proof of Lemma \ref{l3}.
\epf


\subsection{Long time asymptotics}\label{long}

This section is devoted to the study of the asymptotic behavior of system \eqref{e1}, \eqref{de1}
as $t \to \infty$. Put $z_t = \stop_\sigma[s_0,c_1 p + c_2 y]_t = c_1 x_t + c_2 y_t - r_t +\xi_t$.
This allows us, with the notation \eqref{de2}, to rewrite \eqref{e1}, \eqref{de1} in the form
\begin{equation}\label{de3}
\begin{array}{rcl}
(1-a_1)\del_t y + (a_1 + a_2 c_2) y_t + a_2(c_1 - 1) x_t + a_2 s_t  &=& a_2 (z_t-\xi_t) + \epsilon_t,\\
(1-b_1)\del_t x + b_1 s_t - b_2 y_t &=& \eta_t\,.
\end{array}
\end{equation}
As a consequence of \eqref{de3}, we have
$$
(1-b_1)\del_t^2 x + b_1 \del_t s - b_2 \del_t y = \del_t \eta
$$
with the notation \eqref{de2}. This enables us to
eliminate  $y_t$ from the system \eqref{de3} and reformulate it as a second order equation
\be{e4}
\del_t^2 x + A \del_t x + B \del_t s + C x_t + D s_t = h_t\,,\qquad t\ge 2,
\ee
with positive constants
\be{e4b}
A = \frac{a_1 + a_2 c_2}{1-a_1}\,, \ B = \frac{b_1}{1-b_1}\,,
\ C = \frac{a_2 b_2(c_1-1)}{(1-a_1)(1-b_1)}\,, \ D = \frac{b_1(a_1 + a_2 c_2)+ a_2 b_2}{(1-a_1)(1-b_1)}\,,
\ee
and with the right-hand side
\be{deh}
h_t = \frac{1}{(1-a_1)(1-b_1)} \big((1-a_1) \del_t\eta + b_2 (a_2(z_t - \xi_t) + \epsilon_t)
+ (a_1 + a_2 c_2)\eta_t\big).
\ee
The sequence $\{h_t\}$ contains the term $z_t$ which is bounded above by $\sigma$, and the noise terms
$\epsilon_t, \xi_t, \eta_t$, and $\del_t \eta$.

Equation \eqref{e4} always has a solution $x_t$ at each time step $t$, since the right-hand side of \eqref{e4}
is for each fixed $t$ a bounded function of $x_t$ and the left hand side is an increasing piecewise linear
function of $x_t$. In some cases, the solution may not be unique if the coefficient
$D^*:= a_2 b_2/((1-a_1)(1-b_1))$ in front of $x_t$ on the right-hand side is large.
Our computations below show, however, that all solutions have the same asymptotic convergence
towards a small neighborhood of a particular equilibrium point depending on the trajectory.


\subsubsection{Auxiliary estimate 1}\label{est1}

We put
\be{x1}
q_t := C x_t + D s_t\,,
\ee
and multiply the equation \eqref{e4} by $\del_t q = C \del_t x + D \del_t s$.
Putting
\be{de8}
V^1_t := \frac12 \bigr(C (\del_t x)^2 + D (\del_t s)^2 +q_t^2\bigl)
\ee
and using the relations \eqref{de7}--\eqref{de4} we obtain that
\be{de9}
V^1_t - V^1_{t-1}+ \frac{C}2 (\del_t^2 x)^2 + \frac{D}2 (\nabla_t^2 s)^2+ AC (\del_t x)^2 + (BC + AD + BD) (\del_t s)^2
+ \frac12(\del_t q)^2 \le h_t \del_t q\,.
\ee


\subsubsection{Auxiliary estimate 2}\label{est2}

We now rewrite \eqref{e4} in the form
\be{e4a}
\del_t^2 x + \frac{A}{C} \del_t q + \left(B - \frac{AD}{C}\right)\del_t s + q_t = h_t
\ee
with $q_t$ given by \eqref{x1}, and multiply it by $q_t$. We use \eqref{de7} again
and find constants $E, F>0$ depending on $A,B,C,D$ such that
\begin{align} \nonumber
\del_t x\, q_t {-} \del_{t-1} x\, q_{t-1} + \frac{A}{2C} (q_t^2 {-} q_{t-1}^2) + \frac12 q_t^2
&\le (\del_{t} x - \del^2_{t} x) \del_t q + E(\del_t s)^2
- \frac{A}{2C} (\del_t q)^2 + h_tq_t\\ \label{de10}
&\le F \left(\frac12 (\del^2_{t} x)^2 + A(\del_t x)^2\right) + h_t\,q_t\,.
\end{align}
We now set
\be{de11}
V^0_t := \del_t x\, q_t + \frac{A}{2C} q_t^2\,,
\ee
so that \eqref{de10} has the form
\be{de12}
V^0_t - V^0_{t-1} + \frac12 q_t^2 \le F\left(\frac12 (\del^2_{t} x)^2 + A(\del_t x)^2\right) + h_t\,q_t\,.
\ee
Finally, we choose $\lambda>0$ such that $\lambda F < C$, and $\lambda V^0_t \ge -\frac12 V^1_t$
for all $t\in \nat$, and put $W_t = V^1_t + \lambda V^0_t$. Then, due to \eqref{de9} and \eqref{de12}, there exists a constant $\mu>0$
such that
\be{de13}
W_t - W_{t-1} + \mu W_t \le L |h_t|\sqrt{W_t}
\ee
for all $t \in \nat$ as a consequence of \eqref{de8}, \eqref{de9}, and \eqref{de12}, with
some constant $L > 0$.


\subsubsection{Asymptotic behavior: Proof of Theorems \ref{t1}--\ref{t2}}\label{asym}

Let $c_1>1$. If $h_t = 0$, that is, no noise is present and the reaction of the central bank is
instantaneous with $\sigma = 0$, then $W_t$ is a Lyapunov function of the system which decays
exponentially to $0$. In particular, $q_t$ defined by \eqref{x1} converge exponentially to $q_\infty = 0$.
Then, it follows from Lemmas \ref{l2} and \ref{l3} that $x_t$ converge to some value $x_\infty$,
hence $s_t = \stop_\rho[s_0,x]_t$ converge to some $s_\infty$ such that
\be{as1}
Cx_\infty + D s_\infty = 0,
\ee
and $p_t=x_t-s_t\to x_\infty-s_\infty$.
Further, equations \eqref{e1'} imply that $y_t\to y_\infty$ and the point $(x_\infty,y_\infty,p_\infty)$
belongs to the set \eqref{aA}.
This completes the proof of Theorem \ref{t1}.


In the case of a general right-hand side $h_t$, we have for every $t> T >0$
as a consequence of \eqref{de13} that
\be{de14a}
W_t \le (1+\mu)^{-T}W_{t-T} + L \sum_{j=t-T+1}^t |h_j|\sqrt{W_j} (1+\mu)^{j-1-t}\,.
\ee
Assume that there exists $T>0$ such that for all $t_0 \in \nat$ we have
\be{de14c}
\frac{1}{T}\sum_{t=t_0+1}^{t_0+T} |h_j| \le \hat \sigma\,.
\ee
Then it follows from \eqref{de14a} that
\be{de14b}
W_t \le (1+\mu)^{-T}W_{t-T} + \frac{LT\hat\sigma}{1+\mu} \max_{j=t-T+1, \dots, t} \sqrt{W_j}\,
\ee
for all $t>T$. Assume first that there exists $t>T$ such that $W_t \ge W_j$ for all $j=t-T, \dots, t$.
Then \eqref{de14b} yields that
\[
W_t \le (1+\mu)^{-T}W_t + \frac{LT\hat\sigma}{1+\mu} \sqrt{W_t}\,,
\]
hence,
\[
W_t \le \left(\frac{LT \hat\sigma(1+\mu)^{T-1}}{(1+\mu)^T - 1}\right)^2.
\]
This implies in particular that $W_t$ is bounded and we put
\[
W^* := \limsup_{t\to \infty} W_t <\infty.
\]
For an arbitrary $\delta>0$ we find $t_0$ sufficiently large such that for all $t>t_0-T$
we have $W_t \le W^* + \delta$. Then for $t>t_0$ we obtain from \eqref{de14b} that
\[
W_t \le (1+\mu)^{-T}(W^*+\delta) + \frac{LT\hat\sigma}{1+\mu} \sqrt{W^*+\delta}\,,
\]
hence,
\[
W^* \le (1+\mu)^{-T}(W^*+\delta) + \frac{LT\hat\sigma}{1+\mu} \sqrt{W^*+\delta}\,,
\]
that is,
\[
W^*+\delta \le \frac{LT \hat\sigma(1+\mu)^{T-1}}{(1+\mu)^T - 1} \sqrt{W^*+\delta}
+ \frac{\delta(1+\mu)^{T}}{(1+\mu)^T - 1}\,,
\]
and we conclude that
\be{de14i}
W^* \le \left(\frac{LT \hat\sigma(1+\mu)^{T-1}}{(1+\mu)^T - 1}\right)^2 
+ \delta\frac{(1+\mu)^{T}+1}{(1+\mu)^T - 1}\,.
\ee
Since $\delta>0$ is arbitrary, we obtain that
\be{de14}
\limsup_{t\to \infty} W_t \le \left(\frac{LT \hat\sigma(1+\mu)^{T-1}}{(1+\mu)^T - 1}\right)^2.
\ee

Estimate \eqref{de14} gives a uniform upper bound
for the value of $\limsup_{t\to \infty} W_t$ independent
of the initial condition if the noise terms $\epsilon_t, \xi_t, \eta_t$.
In particular, if $|h_t| \le \hat\sigma$ for all $t \in \nat$, then formula \eqref{de14}
holds with $T=1$, that is,
\be{de14j'}
\limsup_{t\to \infty} W_t \le \left(\frac{L \hat\sigma}{\mu}\right)^2 .
\ee
Since by construction $(C x_t+Ds_t)^2=q_t^2\le 2V_t^1\le 4W_t^1$ and by definition $x_*(s_t)=\kappa s_t=-s_tD/C$, it follows that
\be{de14j}
\limsup_{t\to \infty} \,(x_t-x_*(s_t))^2 \le \left(\frac{2L \hat\sigma}{\mu C}\right)^2 .
\ee
Condition \eqref{added} of Theorem \ref{t2} and the definition of $h_t$ (see \eqref{deh}) imply that the estimate  $|h_t| \le \hat\sigma$ holds
with $\hat{\sigma}=L' \sigma + L'' m$, hence \eqref{de14j} implies the first of the estimates \eqref{estdi}.
Further, the first equation in \eqref{e1} implies that
\[
y_t-y_*(s_t)=y_t-\frac{b_1 s_t}{b_2}=\frac{1-b_1}{b_2}\nabla_t x -\frac{\eta_t}{b_2},
\]
hence the second of the estimates \eqref{estdi} follows from the relations $C(\nabla_t x)^2\le 2 V_t^1\le 4 W_t$ and \eqref{de14j'} combined with $\hat{\sigma}= L' \sigma + L'' m$ and \eqref{added}.
This completes the proof of Theorem \ref{t2}.

\medskip
A few remarks are in order. First, estimate \eqref{de14} can be useful if system \eqref{e1} is obtained as  a time discretization of a continuous time system and the increments of noise are of class $\ell^1_{loc}$.
Second, counterparts of Theorems \ref{t1}--\ref{t2} can be obtained by the same method if $p_t$ is allowed to depend also on $y_t$, for example, if
		$p_t = \play_\rho[s_0,x + \delta y]_t$ with a small $\delta>0$.


\subsection{Stability of the multi-agent model}

As mentioned earlier, the proof of Theorem \ref{t3} is parallel to the proof of Theorem \ref{t1} but additionally relies on an inversion formula for the Prandtl-Ishlinskii operator. We start by deriving the inverse operator.

\subsubsection{Inversion of time discrete Prandtl-Ishlinskii operator}\label{praish}


The main tool in our analysis is the following identity which, being inspired by the developments in
\cite[\S 34]{KP}, is essentially due to M. Brokate, see \cite[Proposition 2.2.16]{BS}.



\begin{lemma}[Brokate identity]\label{brok}
	Let $\rho>0$, $\sigma>0$, $\beta \in [-\rho,\rho]$, $\gamma \in [-\sigma,\sigma]$, and $\bbx \in \SS$ be given. 
	Put $\bbxi = \play_\rho[\beta,\bbx]$, $q = \play_\sigma[\gamma,\bbxi]$. Then $q = \play_{\rho+\sigma}[\beta+\gamma,\bbx]$.
	Moreover, the implication
	\be{d3}
	{q}_t - {q}_{t-1} \ne 0 \ \Longrightarrow \ ({p}_t - {p}_{t-1})({q}_t - {q}_{t-1}) > 0
	\ee
	holds for every $t \in \nat$. 
\end{lemma}

\bpf{Proof}\
Put $\vp = \play_{\rho+\sigma}[\beta+\gamma, \bbx]$. We have by definition
\ber\label{d4}
({p}_t -{p}_{t-1}) (x_t - {p}_t - \rho z) &\ge& 0,\\ \label{d5}
({q}_t -{q}_{t-1}) ({p}_t - {q}_t - \sigma z) &\ge& 0,\\ \label{d6}
(\vp_t -\vp_{t-1}) (x_t - \vp_t - (\rho+\sigma) z) &\ge& 0
\eer
for all $t \in \nat$ and all $z \in [-1,1]$. In \eqref{d5}, we may choose $\sigma z = {p}_{t-1} - {q}_{t-1}$
and obtain
$$
({q}_t -{q}_{t-1})^2 \le ({q}_t -{q}_{t-1})({p}_t -{p}_{t-1}),
$$
and the implication \eqref{d3} follows. In particular, the inequality \eqref{d4} remains valid
if we replace ${p}_t -{p}_{t-1}$ with ${q}_t -{q}_{t-1}$, that is,
\be{d7}
({q}_t -{q}_{t-1}) (x_t - {p}_t - \rho z) \ge 0
\ee
for all $t \in \nat$ and all $z \in [-1,1]$. Adding \eqref{d7} to \eqref{d5} yields
\be{d8}
({q}_t -{q}_{t-1}) (x_t - {q}_t - (\rho+\sigma) z) \ge 0
\ee
for all $t \in \nat$ and all $z \in [-1,1]$. We have indeed $|x_t - {q}_t| \le |x_t - {p}_t|+|{p}_t - {q}_t|
\le \rho+\sigma$, hence we may replace $(\rho+\sigma) z$ in \eqref{d6} with $x_t - {q}_t$,
in \eqref{d8} with $x_t - \vp_t$, and sum the two inequalities to obtain
\be{d9}
(({q}_t - \vp_t) - ({q}_{t-1} - \vp_{t-1}))({q}_t - \vp_t) \le 0.
\ee
{}From \eqref{d9} and \eqref{de7} for $z= q-\vp$ it thus follows
$({q}_t - \vp_t)^2 \le ({q}_0 - \vp_0)^2$ for all $t \in \nat$. We have $\vp_0 = x_0 - \beta - \gamma$,
${q}_0 = {p}_0 - \gamma = x_0 - \beta - \gamma = \vp_0$, and this completes the proof. 
\epf

We now recall the definition of  a Prandtl-Ishlinskii operator as a linear combination of play operators.
More specifically, 
let $n \in \nat$, $\nu_0, \dots, \nu_n \in \real$, $0 = \rho_0 < \rho_1 < \dots < \rho_n$,
and $\beta_i \in [-\rho_i, \rho_i]$ be given numbers (in particular, $\beta_0=0$). For $\bbx \in \SS$ and 
$\bbc = (\beta_1, \dots, \beta_n)$ we put
\be{pi1}
\FF[\bbc, \bbx] = \sum_{i=0}^n \nu_i \play_{\rho_i}[\beta_i, \bbx]
\ee
with the convention $\play_0[0,\bbx] = \bbx$. 
The mapping $\FF: [-\rho_1, \rho_1] \times \dots \times [-\rho_n, \rho_n] \times \SS \to \SS$ is called a
(time and memory discrete) {\em Prandtl-Ishlinskii operator\/}.

\begin{hypothesis}\label{h1}
	Let $\FF$ be a Prandtl-Ishlinskii operator as in \eqref{pi1}. The following conditions are assumed to hold:
	\begin{itemize}
		\item[{\rm (i)}] 
		$|\beta_{i} - \beta_{i-1}| \le \rho_{i} - \rho_{i-1}$ for $i = 1, \dots, n$;
		\item[{\rm (ii)}] $A_i:= \sum_{j=0}^i \nu_j > 0$ for $i = 0, \dots, n$.
		\end{itemize}
\end{hypothesis}

We prove the following statement which shows that the inverse of a discrete Prandtl-Ishlinskii operator
is again a Prandtl-Ishlinskii operator. In the continuous case, this result goes back to~\cite{mathz}.
Finite collections of stops have been considered in \cite{equiv}.
A substantially more general situation in the space of regulated functions is considered in \cite{nsf5}.
In fact, the explicit inversion formula presented below can also be deduced from
\cite[Corollary 3.3]{nsf5} which uses deeper results from the Kurzweil integration theory.
In the discrete case, there exists an elementary proof that we present here. 
Note that Lemma \ref{l2} is a special case of Theorem \ref{inver} in the case $n=1$. 

\begin{theorem}\label{inver}
	Let Hypothesis \ref{h1} 
	hold. Define $B_i := 1/A_i$ and
		\[
		\sigma_0 =\gamma_0= 0,\qquad \sigma_{i} - \sigma_{i-1} = A_{i-1} (\rho_{i} - \rho_{i-1}),\qquad
		\gamma_{i} - \gamma_{i-1} =  A_{i-1} (\beta_{i} - \beta_{i-1}),
		\] 
		\[\zeta_0 := B_0,\qquad
	\zeta_i := B_i- B_{i-1}\] 
	for $i=1, \dots, n$. For an arbitrary $\bbx \in \SS$ put
	\[
	\bbv = \sum_{i=0}^n \nu_i \play_{\rho_i}[\beta_i, \bbx].
	\]
	Then, 
	\[
	\bbx = \sum_{i=0}^n \zeta_i \play_{\sigma_i}[\gamma_i, \bbv].
	\]
\end{theorem}

We start with an auxiliary identity.

\begin{proposition}\label{super}
	Let the hypotheses of Theorem \ref{inver} hold. For $i=1, \dots, n$ put
	\[
	\bbx\oi =   \play_{\rho_i - \rho_{i-1}}[\beta_i - \beta_{i-1}, \bbx^{(i-1)}]\,,\ \ \ \ \
	\bbv\oi =   \play_{\sigma_i - \sigma_{i-1}}[\gamma_i - \gamma_{i-1}, \bbv^{(i-1)}]
	\]
	with $\bbx^{(0)} = \bbx$, $\bbv^{(0)} = \bbv$. Then for $j=1, \dots, n$ we have
	\[
	\bbv\oj = A_{j} \bbx\oj + \sum_{i=j+1}^{n} \nu_i \bbx\oi\,.
	\]
\end{proposition}

\bpf{Proof of Proposition \ref{super}}\ 
The definition of the play states that
\begin{align}\label{pr2}
(x\oi_t - x\oi_{t-1}) (x_t^{(i-1)} - x\oi_t - (\rho_{i} - \rho_{i-1})z) &\ge 0\,,\\ \label{pr3}
(v\oi_t - v\oi_{t-1}) (v_t^{(i-1)} - v\oi_t - (\sigma_{i} - \sigma_{i-1})z) &\ge 0
\end{align}
for all $i=1, \dots n$, $t \in \nat$, and $|z| \le 1$. 
For $j=1, \dots, n$ put
\[
(\vp\oj_0, \vp\oj_1, \dots) = \bbphi\oj := A_{j} \bbx\oj + \sum_{i=j+1}^{n} \nu_i \bbx\oi\,.
\]
Then we have
\begin{equation}\label{pr4}
\bbphi^{(0)} = \bbv^{(0)}\,, \quad \vp\oi_t - \vp_t^{(i-1)} = A_{i-1} (x\oi_t - x_t^{(i-1)})\,.
\end{equation}
Hence, we may choose in \eqref{pr2}--\eqref{pr3} suitable values of $z$ in order to obtain
\begin{align}\label{pr5}
(x\oi_t - x\oi_{t-1}) (\vp_t^{(i-1)} - \vp\oi_t - v_t^{(i-1)} + v\oi_t)&\ge 0\,,\\ \label{pr6}
(v\oi_t - v\oi_{t-1}) (v_t^{(i-1)} - v\oi_t - \vp_t^{(i-1)} + \vp\oi_t)&\ge 0
\end{align}
for $i=1, \dots n$ and $t\in \nat$. We now use the implication \eqref{d3} to conclude that for
all $n \ge j \ge i \ge 1$ we have 
\[
\begin{array}{ll}
(x\oj_t - x\oj_{t-1}) (\vp_t^{(i-1)} - \vp\oi_t - v_t^{(i-1)} + v\oi_t)&\ge 0\,,\\ 
(v\oj_t - v\oj_{t-1}) (v_t^{(i-1)} - v\oi_t - \vp_t^{(i-1)} + \vp\oi_t)&\ge 0.
\end{array}
\]
Summing the above inequalities over $i=1, \dots, j$ yields (note that $\bbphi^{(0)} = \bbv^{(0)}$
by \eqref{pr4})
\begin{equation}\label{pr9}
(x\oj_t - x\oj_{t-1}) (v\oj_t - \vp\oj_t) \ge 0, \ \ \ \ \ 
(v\oj_t - v\oj_{t-1}) (\vp\oj_t - v\oj_t) \ge 0
\end{equation}
for $j=1, \dots, n$, $t\in \nat$.
In particular, for $j=n$, we have $\bbphi\on = A_n \bbx\on$. Multiplying the first inequality of \eqref{pr9}
by $A_n$ and adding the second inequality yields
$$
(\vp\on_t - v\on_t - \vp\on_{t-1} + v\on_{t-1}) (\vp\on_t - v\on_t) \le 0
$$
for all $t\in \nat$. Hence, by \eqref{de7},
\be{pr9a}
(\vp\on_t - v\on_t)^2 \le (\vp\on_0 - v\on_0)^2.
\ee
Note that by \eqref{pr4}, $\vp_0\oi - \vp_0^{(i-1)} = A_{i-1} (x\oi_0 - x_0^{(i-1)})
= A_{i-1} (\beta_i - \beta_{i-1}) = \gamma_i - \gamma_{i-1} = v\oi_0 - v_0^{(i-1)}$, hence 
\begin{equation}\label{pr10}
\vp_0\oi = v\oi_0 \quad \for i = 0, \dots n\,,
\end{equation}
and from \eqref{pr9a} it follows that
\begin{equation}\label{pr11}
\bbphi\on = \bbv\on\,.
\end{equation}
We now continue by backward induction and assume that $\bbphi\oj = \bbv\oj$ for some $2\le j \le n$.
{}From \eqref{pr9} written with $j-1$ instead of $j$ we obtain, using the induction hypothesis, that
\begin{align}\label{pr12}
-(x^{(j-1)}_t - x^{(j-1)}_{t-1}) (v\oj_t - \vp\oj_t - v^{(j-1)}_t + \vp^{(j-1)}_t)&\ge 0\,,\\ \label{pr13}
-(v^{(j-1)}_t - v^{(j-1)}_{t-1}) (\vp\oj_t - v\oj_t - \vp^{(j-1)}_t + v^{(j-1)}_t)&\ge 0\,.
\end{align}
We now add \eqref{pr13} to \eqref{pr6}, multiply the sum of \eqref{pr12} with \eqref{pr5}
by $A_{j-1}$, use \eqref{pr4}, and sum the two resulting inequalities to conclude that
$$
\big((v\oj_t - \vp\oj_t - v^{(j-1)}_t + \vp^{(j-1)}_t)
- (v\oj_{t-1} - \vp\oj_{t-1} - v^{(j-1)}_{t-1} + \vp^{(j-1)}_{t-1})\big)
(v\oj_t - \vp\oj_t - v^{(j-1)}_t + \vp^{(j-1)}_t) \le 0.
$$
This is an inequality of type \eqref{de7} which, together with the induction assumption
$\bbphi\oj = \bbv\oj$ yields
\begin{equation}\label{pr14}
(v^{(j-1)}_t - \vp^{(j-1)}_t)^2 \le (v^{(j-1)}_{t-1} - \vp^{(j-1)}_{t-1})^2\,. 
\end{equation}
Referring to \eqref{pr4} and \eqref{pr10}  completes the proof.
\epf

\bpf{Proof of Theorem \ref{inver}}\ 
Let $\bbx\oi, \bbv\oi$ be as in Proposition \ref{super}.
By Lemma \ref{brok} we have $\bbv\oi = \play_{\sigma_i}[\gamma_i, \bbv]$ for all $i = 0,1,\dots, n$.
We now use Proposition \ref{super} and the summation by parts formula, which yields
$$
\sum_{j=0}^{n} \zeta_j \bbv\oi = \sum_{j=0}^{n} \zeta_j \Big(A_{j} \bbx\oj + \sum_{i=j+1}^{n} \nu_i \bbx\oi\Big)
= \bbx^{(0)} + \sum_{i=1}^{n} \bbx\oi (A_{i} \zeta_i +  \nu_i B_{i-1}) = \bbx\,,
$$
and the proof is complete.
\epf


\subsubsection{Proof of Theorem \ref{t3}}\label{inv}
We can now use Theorem \ref{inver} to adapt the proof of Theorem \ref{t1} to the multi-agent case. The 
multi-agent system \eqref{m1}, \eqref{e1'''}
can still be rewritten as the equation
\be{e4p}
\del_t^2 x + A \del_t x + B \del_t s + C x_t + D s_t = 0 
\ee
with positive constants $A, B, C, D$ given by \eqref{e4b}
and  $s_t$ as in \eqref{m3}. From the inequalities \eqref{de4}--\eqref{de6}
we easily obtain their counterparts
\ber \label{m4}
\del_t x \del_t s &=& \nu_0 (\del_t x)^2 + \sum_{i=1}^n \nu_i  \del_t x \del_t s\oi \ge
\nu_0 (\del_t x)^2 + \sum_{i=1}^n \nu_i (\del_t s\oi)^2,\\ \nonumber
\del_t^2 x \del_t s &=& \nu_0 \del_t^2 x \del_t x + \sum_{i=1}^n \nu_i \del_t^2 x \del_t s\oi\ge
\nu_0 \del_t^2 x \del_t x + \sum_{i=1}^n \nu_i \del_t^2 s\oi \del_t s\oi\\ \label{m5}
&\ge& \frac{\nu_0}{2} \big((\del_t x)^2 {-} (\del_{t-1} x)^2 {+} (\del_t^2 x)^2\big)
+ \sum_{i=1}^n\frac{\nu_i}{2} \big((\del_t s\oi)^2 {-} (\del_{t-1} s\oi)^2 {+} (\del_t^2 s\oi)^2\big),\qquad
\eer
where we denote $s\oi := \stop_{\rho_i}[\beta_i, x]$.

As in Section \ref{long}, we define $q_t$ by formula \eqref{x1} and multiply equation \eqref{e4p}
by $\del_t q$. It follows from \eqref{m4}--\eqref{m5} that
\begin{align}\nonumber
\del_t^2 x \del_t q &= (C{+}\nu_0 D) \del_t^2 x \del_t x + D \sum_{i=1}^n \nu_i \del_t^2 x \del_t s\oi\ge
(C{+}\nu_0 D) \del_t^2 x \del_t x + D\sum_{i=1}^n \nu_i \del_t^2 s\oi \del_t s\oi\\ \nonumber 
&\ge \frac{C{+}\nu_0 D}{2} \big((\del_t x)^2 {-} (\del_{t-1} x)^2 {+} (\del_t^2 x)^2\big)
+ \sum_{i=1}^n\frac{\nu_i D}{2} \big((\del_t s\oi)^2 {-} (\del_{t-1} s\oi)^2 {+} (\del_t^2 s\oi)^2\big)
\end{align}
and
\begin{align}\nonumber
(A\del_t x + B \del_t s) \del_t q &\ge (AC{+}\nu_0(BC{+}AD)) (\del_t x)^2+ BD (\del_t s)^2\\ \nonumber 
& + (BC{+}AD)\sum_{i=1}^n \nu_i (\del_t s\oi)^2.
\end{align}
Putting
\[
\tilde V^1_t := \frac12 \left((C{+}\nu_0 D)(\del_t x)^2 + D\sum_{i=1}^n \nu_i (\del_t s\oi)^2 + q_t^2\right),
\]
we obtain (cf.~\eqref{de9} with $h_t=0$)
\begin{align}\nonumber
\tilde V^1_t - \tilde V^1_{t-1} &+ \frac{C{+}\nu_0 D}{2}(\del_t^2 x)^2 +
\sum_{i=1}^n\frac{\nu_i D}{2}(\del_t^2 s\oi)^2 + (AC{+}\nu_0(BC{+}AD)) (\del_t x)^2\\ \nonumber 
&+ BD (\del_t s)^2 + (BC{+}AD)\sum_{i=1}^n \nu_i (\del_t s\oi)^2 + \frac12 (\del_t q)^2 \le 0\,. 
\end{align}
We continue as in Section \ref{long} and multiply \eqref{e4p} by $q_t$. Similarly to \eqref{de10}, we obtain 
\begin{align} \nonumber
\del_t x\, q_t {-} \del_{t-1} x\, q_{t-1} + \frac{A}{2C} (q_t^2 {-} q_{t-1}^2) + \frac12 q_t^2
&\le (\del_{t} x - \del^2_{t} x) \del_t q + E(\del_t s)^2
- \frac{A}{2C} (\del_t q)^2 
\\ \label{m10}
&\le F \left(\frac12 (\del^2_{t} x)^2 + (\del_t x)^2 + \sum_{i=1}^n \nu_i (\del_t s\oi)^2\right) 
\end{align}
with some constants $E, F>0$ depending only on $A, B, C, D$. We set
\[
\tilde V^0_t := \del_t x\, q_t + \frac{A}{2C} q_t^2
\]
and rewrite \eqref{m10} in the form
\[
\tilde V^0_t - \tilde V^0_{t-1} + \frac12 q_t^2
\le F \left(\frac12 (\del^2_{t} x)^2 + (\del_t x)^2 + \sum_{i=1}^n \nu_i (\del_t s\oi)^2\right)
\,,
\]
which is parallel to \eqref{de11}--\eqref{de12}.
We now define an auxiliary energy functional
\[
\tilde W_t = \tilde V^1_t + \lambda \tilde V^0_t
\]
with $\lambda>0$ such that $\lambda \tilde V^0_t \ge -\frac12  \tilde V^1_t$ and
$\lambda < \min\{C+\nu_0 D, AC{+}\nu_0(BC{+}AD), BC{+}AD\}$. We then find $\mu>0$ and $L>0$ such that
for all $t \in \nat$ we have the inequality
%
\[
\tilde W_t - \tilde W_{t-1} + \mu \tilde W_t \le 0.
\]
Therefore, the decay of $\tilde W_t$ is exponential according to the formula
\[
\tilde W_t \le \left(\frac{1}{1+\mu}\right)^t \tilde W_{0},
\]
hence also
\be{m17}
\lim_{t\to \infty} q_t = 0\,.
\ee
Using \eqref{m3}, we can rewrite formula \eqref{x1} as
\[
q_t = (C+ \nu_0 D) x_t + D \sum_{i=1}^n \nu_i \stop_{\rho_i}[\beta_i, x]_t
= (C+D) x_t - D \sum_{i=1}^n \nu_i \play_{\rho_i}[\beta_i, x]_t.
\]
Hence, $q_t$ is given by a Prandtl-Ishlinskii operator of the form
\[
q_t = \tilde \nu_0 x_t + \sum_{i=1}^n \tilde \nu_i \play_{\rho_i}[\beta_i, x]_t
\]
with $\tilde \nu_0 = C+D$, $\tilde \nu_i = - \nu_i D$ for $i = 1, \dots, n$. The hypotheses of
Theorem \ref{inver} are satisfied for
$$
\tilde A_i := \sum_{j=0}^i \tilde \nu_j \ge C + D\left(1 - \sum_{j=0}^i \nu_j\right) \ge C > 0.
$$
Consequently, by Theorem \ref{inver}, we have
\[
x_t = \tilde \zeta_0 q_t + \sum_{i=1}^n \tilde \zeta_i \play_{\sigma_i}[\gamma_i, q]_t
\]
with suitable constants $\tilde \zeta_i, \sigma_i, \gamma_i$. From \eqref{m17} and from Lemma \ref{l3} we
conclude that $x_t$ and $s_t$ are Cauchy sequences that converge to some limits $x_\infty$, $s_\infty$,
respectively, which completes the proof of Theorem \ref{t3}.


\section{Conclusions}

We have replaced rational expectations about future inflation 
with 
a form of
boundedly rational aggregated `sticky' expectation modeled by the play operator in a simple standard macroeconomic
model. This single (and conceptually quite
elementary) change transforms a unique equilibrium linear system to a PWL system with an entire
continuum of 
equilibrium states. The PWL model with $n$ agents has $2n$ switching surfaces and an $n$-dimensional continuum of equilibria.
By constructing a Lyapunov function and developing a technique for inverting the Prandtl-Ishlinskii operator, we have shown that, when there is no exogenous noise, the continuum of equilibria is the global attractor of the system.
The size of the basin of attraction of a particular equilibrium varies, generally
becoming smaller towards the boundary of the set of equilibrium states.

If the presence of stickiness/frictions in economics does indeed
induce a myriad of coexisting (metastable) equilibria then phenomena that are not
possible (or require a posteriori model adjustments) in unique
equilibrium models become not just feasible but inevitable. Perhaps
the most obvious of these permanence, also known as
remanence, where a system does not revert to its previous state after
an exogenous shock is applied and then removed. It is of course a central concern
of macreconomics whether or not economies affected by, say, significant
negative shocks can be expected to have permanently reduced
productivity levels.
For the models studied in this paper, 
after sufficiently small shocks (whether exogenous or applied by policy makers) the system
will indeed revert to the same equilibrium but larger shocks will move the system from the basin of attraction of one equilibrium to the basin of attraction of a different one (at the same model parameters).  The path to this new equilibrium may be long with a highly unpredictable endpoint.
Furthermore, in
the latter case 
the system 
will not
exhibit a tendency to
return to its pre-shocked state --- the model displays true permanence. And 
the model parameters alone cannot determine which
equilibrium a system is currently in without knowing important
information about the prior states of the system --- true path dependence. 
Hence, the model accounts for several hard-to-explain empirical regularities observed in economic data.
This feature of the model is significant not just because it
 corresponds closely to actual economic events but it may have
 implications for forecasting and policy prescriptions too. 

Our model of expectation formation is thus both mathematically
tractable and has some basis in both observed data (see also \cite{g, laura})
and models of bounded
rationality. As such it provides a potentially useful,
analytically tractable, alternative
to staggered/delayed models --- and one with additional
complexity and explanatory power. 
Our choice of inflation expectations as the candidate for an initial
investigation was influenced by the work of De Grauwe \cite{grauwe} on
a different type of boundedly rational expectation formation process
in a simple DSGE model.  However, play operators are also a viable
candidate for modeling other sticky economic variables at both the
micro- and macro-economic levels. To demonstrate this, 
we used a play operator to represent sticky responses by the Central Bank.
Although it has not been
relevant to this paper play and stop operators, when combined appropriately
\cite{nsf5} can have a remarkably simple aggregated response, even when
connected via a network. This allows for (almost)-analytic solutions
even when cascades and rapid transitions between states are occurring
and will be the subject of future work.
The same form of stickiness described above with the associated 
	play operators have already been used to develop
	non-equilibrium asset-pricing models \cite{nsf1}.
		
\section*{Acknowledgments}
DR was supported by NSF grant DMS-1413223, PK was supported by the GA\v CR Grant GA15-12227S
and RVO: 67985840. The hospitality of the University of Texas at Dallas during PK's stay there
in November 2015 is gratefully acknowledged.

\bibliographystyle{elsarticle-num}

\bibliography{kkpr_bibJun26}


\end{document}